\documentclass[12pt,a4paper]{article}
\usepackage[cp1251]{inputenc} 
\usepackage{amsfonts}

\newcommand{\di}{\displaystyle}

\newcommand{\al}{\alpha}
\newcommand{\be}{\beta}
\newcommand{\ga}{\gamma}
\newcommand{\de}{\delta}
\newcommand{\la}{\lambda}
\newcommand{\om}{\omega}

\newcommand{\vv}{\varphi}
\newcommand{\iy}{\infty}

\begin{document}

\begin{center}
{\large\bf
Inverse Spectral Problems for Sturm-Liouville Operators on Hedgehog-type Graphs with General Matching Conditions}\\[0.2cm]
{\bf V.\,Yurko} \\[0.2cm]
\end{center}

\thispagestyle{empty}

\noindent {\bf Abstract.} Boundary value problems on hedgehog-type graphs
for Sturm-Liouville differential operators with general matching conditions are studied.
We investigate inverse spectral problems of recovering the coefficients of the differential 
equation from the spectral data. For this inverse problem we prove a uniqueness theorem 
and provide a procedure for constructing its solution.



\noindent {\bf 1. Introduction. }\\
{\bf 1.1. } We study an inverse spectral problem for Sturm-Liouville differential operators 
on so-called hedgehog-type graphs with general matching conditions in the interior
vertices. Inverse spectral problems consist in recovering operators from their spectral
characteristics. The main results on inverse spectral problems for Sturm-Liouville
operators on an {\it interval} are presented in the monograph [1] and other works.
Differential operators on graphs (networks, trees) often appear in natural sciences
and engineering (see [2] and the references therein). Most of the results in this
direction are devoted to direct problems of studying properties of the
spectrum and the root functions for operators on graphs. Inverse spectral problems,
because of their nonlinearity, are more difficult to investigate, and nowadays there
exists only a small number of papers in this area. In particular, inverse spectral problems
of recovering {\it coefficients} of differential operators {\it on trees} (i.e on
graphs without cycles) were solved in [3-4]. Inverse problems for Sturm-Liouville
operators on graphs with a cycle were studied in [5] and other papers but only
in the case of so-called {\it standard matching conditions.} 

In the present paper we consider Sturm-Liouville operators on hedgehog-type graphs 
with generalized matching conditions. This class of matching 
conditions appears in applications and produces new qualitative difficulties in 
investigating nonlinear inverse coefficient problems. For studying this class of inverse 
problems we develop the ideas of the method of spectral mappings [6]. We prove 
a uniqueness theorem for this class of nonlinear inverse problems, and provide
a constructive procedure for their solution. In order to construct the solution, we solve, 
in particular, an important auxiliary inverse problem for a quasi-periodic boundary value
problem on the cycle with discontinuity conditions in interior points. The obtained results 
are natural generalizations of the well-known results on inverse problems for differential 
operators on an interval and on graphs with standard matching conditions.

\medskip
{\bf 1.2. } Consider a compact graph $G$ in ${\bf R^m}$ with the set of edges
${\cal E}=\{e_0,\ldots, e_r\},$ where $e_0$ is a cycle, ${\cal E}'=\{e_1,\ldots, e_r\}$
are boundary edges. Let $\{v_1,\ldots, v_{r+N}\}$ be the set of vertices,
where $V=\{v_1,\ldots, v_r\},$ $v_k\in e_k$, are boundary vertices, and
$U=\{v_{r+1},\ldots, v_{r+N}\}$ are internal vertices, $U={\cal E}'\cap e_0$.
The cycle $e_0$ consists of $N$ parts:
$e_0=\bigcup_{k=1}^N e_{r+k}$,  $e_{r+k}=[v_{r+k},v_{r+k+1}],\;
k=\overline{1,N},\; v_{r+N+1}:=v_{r+1}.$
Each boundary edge $e_j$, $j=\overline{1,r}$ has the initial point at $v_j$,
and the end point in the set $U.$ The set ${\cal E}'$ consists of $N$ groups
of edges:
${\cal E}'={\cal E}_1\cup\ldots\cup {\cal E}_N\,, {\cal E}_k\cap e_0=v_{r+k}.$
Let $r_k$ be the number of edges in ${\cal E}_k$; hence $r=r_1+\cdots+r_N$.
Denote $m_0=1,$ $m_k=r_1+\cdots+r_k$, $k=\overline{1,N}.$ Then
${\cal E}_k=\{e_j\},\; j=\overline{m_{k-1}+1, m_k},$
$v_{r+k}=\bigcap_{j=m_{k-1}+1}^{m_k} e_j\,,\; k=\overline{1,N}.$
Thus, the boundary edge $e_j\in {\cal E}_k$ can be viewed as the segment
$e_j=[v_j, v_{r+k}].$ 

Let $T_j$ be the length of the edge $e_j$, $j=\overline{1,r+N},$ and let
$T:=T_{r+1}+\ldots +T_{r+N}$ be the length of the cycle $e_0$. Put
$b_0=0,\; b_k=T_{r+1}+\ldots +T_{r+k},$ $k=\overline{1,N}.$ Then $b_N=T.$

Each edge $e_j$, $j=\overline{1,r+N}$ is parameterized by the parameter
$x_j\in [0,T_j],$ and $x_j=0$ corresponds to the vertex $v_j$.
The whole cycle $e_0$ is parameterized by the parameter $x\in [0,T],$ where
$x=x_{r+j}+b_{j-1}$ for $x_{r+j}\in[0, T_{r+j}],$ $j=\overline{1,N}.$

An integrable function $Y$ on $G$ may be represented as $Y=\{y_j\}_{j=\overline{1,r+N}}$,
where the function $y_j(x_j),$ $x_j\in [0,T_j],$ is defined on the edge $e_j$.
The function $y(x),\; x\in[0,T]$ on the cycle $e_0$ is defined by
$y(x)=y_{r+j}(x_{r+j}),$ $j=\overline{1,N}.$

Let $Q=\{q_j\}_{j=\overline{1,r+N}}$ be an integrable real-valued
function on $G$; $Q$ is called the potential. The function $q(x),\; x\in[0,T]$ is 
defined by $q(x)=q_{r+j}(x_{r+j}),$ $j=\overline{1,N}.$ Denote 
$U_j(Y):=y'_j(0)-h_jy_j(0),$ $j=\overline{1,r+N},$ $U_{r+N+1}:=U_{r+1},$
where $h_j$ are real numbers. Consider the following differential equation on $G$:
$$
-y''_j(x_j)+q_j(x_j)y_j(x_j)=\la y_j(x_j),\quad x_j\in [0,T_j],\quad j=\overline{1,r+N},  \eqno(1)
$$
the functions $y_j, y'_j,$ $j=\overline{1,r+N},$ are absolutely continuous on $[0,T_j]$ 
and satisfy the following matching conditions in each internal vertex $v_{\mu+1}$, 
$\mu=\overline{r+1,r+N}$:
$$
y_{\mu+1}(0)=\al_j y_{j}(T_{j})\;\mbox{for all}\; e_j\in {\cal E}'_{\mu-r+1},\quad
U_{\mu+1}(Y)=\sum_{e_j\in {\cal E}'_{\mu-r+1}} \be_j y'_j(T_j),                                  \eqno(2)
$$
$$
y_{r+N+1}:=y_{r+1},\; h_{r+N+1}:=h_{r+1},\; {\cal E}_{N+1}:={\cal E}_1,\;
{\cal E}'_{\mu-r+1}:={\cal E}_{\mu-r+1}\cup e_{\mu},
$$
where $\al_j$ and $b_j$ are real numbers, and $\al_j\be_j\ne 0.$ For definiteness, let
$\al_j\be_j>0.$ The matching conditions (2) are a generalization of the standard matching
conditions (see [5]), where $\al_j=\be_j=1,$ $h_j=0.$

Let us consider the boundary value problem $B_0$ on $G$ for equation (1) with
the matching conditions (2) and with the following boundary conditions at the
boundary vertices $v_1,\ldots, v_r$:
$$
U_j(Y)=0,\quad j=\overline{1,r}.
$$
Denote by $\Lambda_0=\{\la_{n0}\}_{n\ge 0}$ the eigenvalues (counting with multiplicities)
of $B_0$. Moreover, we also consider the boundary value problems $B_{\nu_1,\ldots,\nu_p},$
$p=\overline{1,r},$ $1\le\nu_1<\ldots\nu_p\le r$ for equation (1) with the matching
conditions (2) and with the boundary conditions
$$
y_k(0)=0,\; k=\nu_1,\ldots,\nu_p, \quad
U_j(Y)=0,\; j=\overline{1,r},\; j\ne\nu_1,\ldots,\nu_p.
$$
Denote by $\Lambda_{\nu_1,\ldots,\nu_p}:=\{\la_{n,\nu_1,\ldots,\nu_p}\}_{n\ge 0}$
the eigenvalues (counting with multiplicities) of $B_{\nu_1,\ldots,\nu_p}$.

It will be shown in Section 2 that an important role for solving inverse problems
on graphs is played by an auxiliary quasi-periodic boundary value problem on the cycle
with discontinuity conditions in interior points. The parameters of this auxiliary
problem depend on the parameters of $B_0$. More precisely, let us introduce
real numbers $\ga_j, \eta_j,\;(j=\overline{1,N-1}),\;h, \al, \be:$ 
$$
\left.\begin{array}{c}
\ga_j=\di\sqrt{\di\frac{\al_{r+j}}{\be_{r+j}}}\,,\; \eta_j=\ga_j h_{r+j+1},
\; j=\overline{1,N-1},\quad h=h_{r+1},\\[4mm]
\al=\al_{r+N}\di\prod_{j=1}^{N-1} \ga_j \di\prod_{j=1}^{N-1} \be_{r+j},
\quad \be=\di\prod_{j=1}^{N-1} \ga_j \di\prod_{j=1}^{N} \be_{r+j}.
\end{array}\right\}                                                                                                        \eqno(3)
$$
Clearly, $\al\be>0,\; \ga_j>0,\; j=\overline{1,N-1}.$ Using these parameters
we consider the following quasi-periodic discontinuity boundary value problem 
$B$ on the cycle $e_0$:
$$
-y''(x)+q(x)y(x)=\la y(x),\quad x\in(0,T),                                                                \eqno(4)
$$
$$
y(0)=\al y(T),\quad y'(0)-hy(0)=\be y'(T),                                                                \eqno(5)
$$
$$
y(b_j+0)=\ga_jy(b_j-0),\;
y'(b_j+0)=\ga_j^{-1}y'(b_j-0)+\eta_jy(b_j-0),\;j=\overline{1,N-1},                    \eqno(6)
$$
$$
0<b_1<\ldots <b_{N-1}<b_N=T.
$$
Let $S(x,\la)$ and $C(x,\la)$ be solutions of equation (4) satisfying discontinuity
conditions (6) and the initial conditions $S(0,\la)=C'(0,\la)=0,$ $S'(0,\la)=C(0,\la)
=1.$ Put $\vv(x,\la)=C(x,\la)+hS(x,\la).$ Eigenvalues $\{\la_n\}_{n\ge 1}$ of 
$B$ coincide with zeros of the characteristic function
$$
a(\la)=\al\vv(T,\la)+\be S'(T,\la)-(1+\al\be).                                                      \eqno(7)
$$
Put $d(\la):=S(T,\la),$ $Q(\la)=\al\vv(T,\la)-\be S'(T,\la).$
All zeros $\{z_n\}_{n\ge 1}$ of the entire function $d(\la)$ are simple, 
i.e. $\dot d(z_n)\ne 0,$ where $\dot d(\la):=\frac{d}{d\la}\,d(\la).$
Denote $M_n=-d_1(z_n)/\dot d(z_n),$ where $d_1(\la):=C(T,\la).$
The sequence $\{M_n\}_{n\ge 1}$ is called the Weyl sequence. Put
$\omega_n=\mbox{sign}\,Q(z_n),\; \Omega=\{\omega_n\}_{n\ge 1}.$

We choose and fix one edge $e_{\xi_i}\in {\cal E}_i$ from each
block ${\cal E}_i$, $i=\overline{1,N},$ i.e. $m_{i-1}+1\le\xi_i\le
m_i$. Denote by $\xi:=\{k:\; k=\xi_1,\ldots,\xi_N\}$ the set of
indices $\xi_i$, $i=\overline{1,N}.$ Let $\al_j$ and $\be_j$,
$j=\overline{1,r+N},$ be known a priori. The inverse problem is
formulated as follows.

\smallskip
{\bf Inverse problem 1.} Given $2^N+r-N$ spectra $\Lambda_j,\;
j=\overline{0,r}$, $\Lambda_{{\nu_1,\ldots,\nu_p}}$,
$p=\overline{2,N},$ $1\le\nu_1<\ldots<\nu_p\le r,$ $\nu_j\in\xi,$
and $\Omega,$ construct the potential $Q$ on $G$ and 
$H:=[h_j]_{j=\overline{1,r+N}}$.

\smallskip
Obviously, in general it is not possible to recover also all coefficients $\al_j$
and $\be_j$. Note that this inverse problem is a generalization of the classical
inverse problems for Sturm-Liouville operators on an interval or on graphs.

\smallskip
Let us formulate the uniqueness theorem for the solution of
Inverse Problem 1. For this purpose together with $q$ we consider
a potential $\tilde q.$ Everywhere below if a symbol $\al$ denotes
an object related to $q,$ then $\tilde\al$ will denote the
analogous object related to $\tilde q.$

\smallskip
{\bf Theorem 1. }{\it If $\Lambda_j=\tilde\Lambda_j,\; j=\overline{0,r},$
$\Lambda_{{\nu_1,\ldots,\nu_p}}=\tilde\Lambda_{{\nu_1,\ldots,\nu_p}}$,
$p=\overline{2,N},$ $1\le\nu_1<\ldots<\nu_p\le r,$ $\nu_j\in\xi,$
and $\Omega=\tilde\Omega,$ then $Q=\tilde Q$ and $H=\tilde H.$}

\smallskip
We will also provide there a constructive procedure for the solution of 
Inverse problem 1. 

\medskip
\noindent {\bf 2. Solution of the inverse problem.}\\
{\bf 2.1. } Let $S_j(x_j,\la),\;C_j(x_j,\la),\;j=\overline{1,r+N},\;x_j\in[0,T_j],$
be the solutions of equation (1) on the edge $e_j$ with the initial conditions
$S_j(0,\la)=C'_j(0,\la)=0,$ $S'_j(0,\la)=C_j(0,\la)=1.$
Put $\vv_j(x,\la)=C_j(x_j,\la)+h_jS_j(x_j,\la).$ Clearly,
$\langle \vv_j(x_j,\la),\, S_j(x_j,\la)\rangle\equiv 1,$
where $\langle y,z\rangle :=yz'-y'z$ is the Wronskian of $y$ and $z.$

Let here and below $\la=\rho^2,\; \tau:=Im\,\rho\ge 0,$
$\Pi:=\{\rho:\; \tau\ge 0\},$ $\Pi_{\de}:=\{\rho:\;
\arg\rho\in[\de,\pi-\de]\},$ $\de\in (0,\pi/2).$ 
For $|\rho|\to\iy,\; \rho\in\Pi_{\de}$ the following relations hold (see [7]):
$$
a(\la)=\frac{(\al+\be)\xi}{2}\,e^{-i\rho T}[1],\quad d(\la)=
-\frac{\xi}{2i\rho}\,e^{-i\rho T}[1],\quad \xi:=\prod_{j=1}^{N-1} \xi_j^{+}\,,     \eqno(8)
$$
$$
a(\la)=O(e^{\tau T}),\quad
d(\la)=O(\rho^{-1}e^{\tau T}),\quad |\rho|\to\iy,\quad \rho\in \Pi.                      \eqno(9)
$$

Fix $k=\overline{1,r}.$ Let $\Phi_k=\{\Phi_{kj}\}_{j=\overline{1,r+N}}$,
be the solution of equation (1) satisfying (2) and the boundary conditions
$$
U_j(\Phi_{k})=\de_{jk}, \quad j=\overline{1,r},                                              \eqno(10)
$$
where $\de_{jk}$ is the Kronecker symbol. Denote $M_k(\la):=\Phi_{kk}(0,\la),$
$k=\overline{1,r}.$ The function $M_k(\la)$ is called the {\it Weyl function}
with respect to the boundary vertex $v_k$. Clearly,
$$
\Phi_{kk}(x_k,\la)=S_{k}(x_k,\la)+M_k(\la)\vv_{k}(x_k,\la),
\quad x_k\in[0,T_k],\quad k=\overline{1,r},                                               \eqno(11)
$$
$$
\langle\vv_{k}(x_k,\la),\,\Phi_{kk}(x_k,\la)\rangle\equiv 1.                      \eqno(12)
$$
Denote $M_{kj}^1(\la):=\Phi_{kj}(0,\la)$, $M_{kj}^0(\la):=\Phi'_{kj}(0,\la)$. Then
$$
\Phi_{kj}(x_j,\la)=M_{kj}^0(\la)S_j(x_j,\la)+M_{kj}^1(\la)\vv_j(x_j,\la),
\quad x_j\in[0,T_j],\; j=\overline{1,r+N},\; k=\overline{1,r}.                    \eqno(13)
$$
In particular, $M_{kk}^0(\la)=1$, $M_{kk}^1(\la)=M_k(\la).$ Substituting
(13) into (2) and (10) we obtain a linear algebraic system $D_k$ with respect
to $M_{kj}^\nu(\la)$, $\nu=0,1,\; j=\overline{1,r+N}.$ The determinant
$\Delta_0(\la)$ of $D_k$ does not depend on $k$ and has the form
$$
\Delta_0(\la)=\sigma(\la)\Big(a_0(\la)+\di\sum_{k=1}^{N}
\di\sum_{1\le\mu_1<\ldots<\mu_k\le N} a_{\mu_1\ldots\mu_k}(\la)
\prod_{i=1}^k\Big(\sum_{e_j\in{\cal E}_{\mu_i}}\Omega_j(\la)\Big)\Big),         \eqno(14)
$$
$$
\sigma(\la)=\prod_{j=1}^r (\al_j\vv_j(T_j,\la)),\quad
\Omega_j(\la)=\frac{\be_j\vv'_j(T_j,\la)}{\al_j\vv_j(T_j,\la)}\,,  
\quad a_0(\la)=a(\la),\quad a_1(\la)=\al d(\la).                                                      \eqno(15)
$$
We note that the coefficients $a_0(\la)$ and $a_{\mu_1\ldots\mu_k}(\la)$
in (14) depend only on $S_j^{(\nu)}(T_j,\la)$ and $C_j^{(\nu)}(T_j,\la),$
for $j=\overline{r+1,r+N}.$ We do not need
concrete formulae for the other coefficients $a_{\mu_1\ldots\mu_k}(\la).$
The function $\Delta_0(\la)$ is entire in $\la$ of order $1/2$, and its
zeros coincide with the eigenvalues of the boundary value problem $B_0.$
The function $\Delta_0(\la)$ is called the characteristic function for
the boundary value problems $B_0.$ Let $\Delta_{\nu_1,\ldots,\nu_p}(\la),$
$p=\overline{1,r},$ $1\le\nu_1<\ldots<\nu_p\le r,$ be the function obtained
from $\Delta_0(\la)$ by the replacement of $\vv_j^{(\nu)}(T_j,\la)$ with
$S_j^{(\nu)}(T_j,\la)$ for $j=\nu_1,\ldots,\nu_p$, $\nu=0,1.$ More precisely,
$$
\Delta_{\nu_1,\ldots,\nu_p}(\la)=\sigma_{\nu_1,\ldots,\nu_p}(\la)\Big(a_0(\la)+
\di\sum_{k=1}^{N} \di\sum_{1\le\mu_1<\ldots<\mu_k\le N} a_{\mu_1\ldots\mu_k}(\la)
$$
$$
\times\prod_{i=1}^k\Big(\sum_{e_j\in{\cal E}_{\mu_i},\;j\ne\nu_1,\ldots,\nu_p}
\Omega_j(\la)+\sum_{e_j\in{\cal E}_{\mu_i},\;
j=\nu_1,\ldots,\nu_p} \Omega_j^0(\la)\Big)\Big),                                             \eqno(16)
$$
$$
\sigma_{\nu_1,\ldots,\nu_p}(\la)=\prod_{j=1,\; j\ne \nu_1,\ldots,\nu_p}^r
(\al_j\vv_j(T_j,\la)) \prod_{j=\nu_1,\ldots,\nu_p} (\al_j S_j(T_j,\la)),\;
\Omega_j^0(\la)=\frac{\be_j S'_j(T_j,\la)}{\al_j S_j(T_j,\la)}\,.                \eqno(17)
$$
The function $\Delta_{\nu_1,\ldots,\nu_p}(\la)$ is entire in $\la$ of order
$1/2$, and its zeros coincide with the eigenvalues of the boundary value problem
$B_{\nu_1,\ldots,\nu_p}.$ The function $\Delta_{\nu_1,\ldots,\nu_p}(\la)$ is called
the characteristic function for the boundary value problem $B_{\nu_1,\ldots,\nu_p}.$

Solving the algebraic system $D_k$ we get by Cramer's rule: $M_{kj}^s(\la)
=\Delta_{kj}^s(\la)/\Delta_0(\la)$, $s=0,1,\,j=\overline{1,r+N},$ where
the determinant $\Delta_{kj}^s(\la)$ is obtained from $\Delta_0(\la)$ by
the replacement of the column which corresponds to $M_{kj}^s(\la)$ with
the column of free terms. In particular,
$$
M_k(\la)=-\Delta_k(\la)/\Delta_0(\la),\quad k=\overline{1,r}.        \eqno(18)
$$

{\bf 2.2. } It is known (see [8]) that for each fixed $j=\overline{1,r+N},$
on the edge $e_j$, there exists a fundamental system of solutions of
equation (1) $\{e_{j1}(x_j,\rho), e_{j2}(x_j,\rho)\},$ $x_j\in [0,T_j],
\,\rho\in\Pi,\, |\rho|\ge\rho^*$ with the properties:\\
1) the functions $e_{js}^{(\nu)}(x_j,\rho),\,\nu=0,1,$ are continuous
for $x_j\in [0,T_j],\, \rho\in\Pi,\, |\rho|\ge\rho^*$;\\
2) for each $x_j\in[0,T_j],$ the functions $e_{js}^{(\nu)}(x_j,\rho),\,
\nu=0,1,$ are analytic for $\mbox{Im}\,\rho>0,\,|\rho|>\rho^*$;\\
3) uniformly in $x_j\in[0,T_j],$ the following asymptotical formulae hold
$$
e_{j1}^{(\nu)}(x_j,\rho)=(i\rho)^{\nu}\exp(i\rho x_j)[1],\;
e_{j2}^{(\nu)}(x_j,\rho)=(-i\rho)^{\nu}\exp(-i\rho x_j)[1],
\; \rho\in\Pi,\; |\rho|\to\iy,                                                                    \eqno(19)
$$
where $[1]=1+O(\rho^{-1}).$ Fix $k=\overline{1,r}.$ One has
$$
\Phi_{kj}(x_j,\la)=A_{kj}^{1}(\rho)e_{j1}(x_j,\rho)+
A_{kj}^{0}(\rho)e_{j2}(x_j,\rho),\quad x_j\in[0,T_j],\quad j=\overline{1,r+N}. \eqno(20)
$$
Substituting (20) into (2) and (10) we obtain a linear algebraic system
$D_k^0$ with respect to $A_{kj}^\nu(\rho)$, $\nu=0,1,\; j=\overline{1,r+N}.$
The determinant $\delta(\rho)$ of $D_k^0$ does not depend on $k,$ and has
the form
$$
\de(\rho)=\Big(\de_0+O\Big(\frac{1}{\rho}\Big)\Big)\rho^{r+N}
\exp\Big(-i\rho \sum_{j=1}^{r+N} T_j\Big),                                           \eqno(21)
$$
where $\de_0$ is the determinant obtained from $\de(\rho)$ by the
replacement of $e_{j1}^{(\nu)}(0,\rho),$
$e_{j1}^{(\nu)}(T_j,\rho),$ $e_{j2}^{(\nu)}(0,\rho),$
$e_{j2}^{(\nu)}(T_j,\rho)$ and $h_j$ with $1, 0, (-1)^\nu,
(-1)^\nu$ and $0,$ respectively. We assume that $\de_0\ne 0.$ This
condition is called the {\em regularity condition} for matching.
Differential operators on $G$ which do not satisfy the regularity
condition, possess qualitatively different properties in
connection with the formulation and investigation of inverse
problems, and are not considered in this paper; they require a
separate investigation. We note that for classical Kirchhoff's
matching conditions we have $\al_{j}=\be_{j}=1,\, h_{j}=0,$ and
the regularity condition is satisfied obviously. Solving the
algebraic system $D_k^0$ and using (19)-(21) we get for each fixed
$x_k\in[0,T_k)$:
$$
\Phi_{kk}^{(\nu)}(x_k,\la)=(i\rho)^{\nu-1}\exp(i\rho x_k)[1],
\quad \rho\in\Pi_\de,\; |\rho|\to\infty.                                                    \eqno(22)
$$
In particular, $M_k(\la)=(i\rho)^{-1}[1],\;\rho\in\Pi_\de,
\; |\rho|\to\infty.$ Moreover, uniformly in $x_k\in[0,T_k],$
$$
\vv_k^{(\nu)}(x_k,\la)
=\frac{1}{2}\Big((i\rho)^{\nu}\exp(i\rho x_k)[1]
+(-i\rho)^{\nu}\exp(-i\rho x_k)[1]\Big),\;\rho\in\Pi,\;|\rho|\to\iy.         \eqno(23)
$$
Using (14), (23), (8) and (9), by the well-known method, 
one can obtain the following properties of
the characteristic function $\Delta_0(\la)$ and the eigenvalues
$\Lambda_0$ of the boundary value problem $B_0$.

1) For $\rho\in\Pi,\;|\rho|\to\infty$:
$\Delta_0(\lambda)=O(\exp(\tau \sum_{j=1}^{r+N} T_j)).$

2) There exist $h>0,\; C_h>0$ such that
$|\Delta_0(\lambda)|\ge C_h\exp(\tau \sum_{j=1}^{r+N} T_j)$
for $\tau\ge h.$ Hence, the eigenvalues $\lambda_{n0}=\rho_{n0}^2$
lie in the domain $0\le\tau<h.$

3) The number $N_{\xi}$ of zeros of $\Delta_0(\lambda)$ in the
rectangle $\Lambda_\xi = \{\rho:\; \tau\in[0,h],\; \mbox{Re}\,
\rho\in[\xi,\xi+1]\}$ is bounded with respect to $\xi.$

4) For $n\to\infty,$
$\rho_{n0}=\rho_{n0}^0 + O(1/\rho_{n0}^0),$
where $\la_{n0}^0=(\rho_{n0}^0)^2$ are the eigenvalues of the
boundary value problem $B_0$ with $Q=0$ and $H=0.$

The characteristic functions $\Delta_{\nu_1,\ldots,\nu_p}(\la)$ have
similar properties. In particular, for $\rho\in\Pi,\;|\rho|\to\infty,$
$\Delta_{\nu_1,\ldots,\nu_p}(\lambda)=O(|\rho|^{-p}
\exp(\tau \sum_{j=1}^{r+N} T_j)).$

Using the properties of the characteristic functions and Hadamard's
factorization theorem, one gets that the specification of
the spectrum $\Lambda_0$ uniquely determines the characteristic function
$\Delta_0(\la)$, i.e. if $\Lambda_0=\tilde\Lambda_0$, then
$\Delta_0(\la)\equiv\tilde\Delta_0(\la).$ Analogously, if
$\Lambda_{\nu_1,\ldots,\nu_p}=\tilde\Lambda_{\nu_1,\ldots,\nu_p}$, then
$\Delta_{\nu_1,\ldots,\nu_p}(\la)\equiv\tilde\Delta_{\nu_1,\ldots,\nu_p}(\la).$
The characteristic functions can be constructed as the corresponding
infinite products (see [1] for details).

\medskip
{\bf 2.3. } Fix $k=\overline{1,r},$ and consider the following auxiliary
inverse problem on the edge $e_k$, which is called IP(k).

\smallskip
{\bf IP(k).} Given two spectra $\Lambda_0$ and $\Lambda_k$, construct
$q_k(x_k),\, x_k\in[0,T_k],$ and $h_{k}$.

\smallskip
In IP(k) we construct the potential only on the edge $e_k$, but the spectra
bring a global information from the whole graph. In other words, IP(k) is
not a local inverse problem related to the edge $e_k$.
Let us prove the uniqueness theorem for the solution of IP(k).

\smallskip
{\bf Theorem 2. }{\it Fix $k=\overline{1,r}.$ If $\Lambda_0=\tilde\Lambda_0$ and
$\Lambda_k=\tilde\Lambda_k$, then $q_k(x_k)=\tilde q_k(x_k)$, a.e. on $[0,T_k],$
and $h_{k}=\tilde h_{k}$. Thus, the specification of two spectra $\Lambda_0$ and
$\Lambda_k$ uniquely determines the potential $q_k$ on the edge $e_k$, and the
coefficient $h_{k}$.}

\smallskip
{\it Proof.}  Since $\Lambda_0=\tilde\Lambda_0,\, \Lambda_k=\tilde\Lambda_k$,
it follows that
$\Delta_0(\la)\equiv\tilde\Delta_0(\la),\;\Delta_k(\la)\equiv\tilde\Delta_k(\la),$
and according to (18),
$$
M_k(\la)=\tilde M_k(\la).                                                                  \eqno(24)
$$
Consider the functions
$$
P^k_{1s}(x_k,\la)=(-1)^{s-1}\Big(\vv_{k}(x_k,\la)\tilde\Phi^{(2-s)}_{kk}(x_k,\la)
-\tilde\vv^{(2-s)}_{k}(x_k,\la)\Phi_{kk}(x_k,\la)\Big),\quad s=1,2.           \eqno(25)
$$
Using (12) we calculate
$$
\vv_k(x_k,\la)=P^k_{11}(x_k,\la)\tilde \vv_k(x_k,\la)+
P^k_{12}(x_k,\la)\tilde\vv'_k(x_k,\la).                                       \eqno(26)
$$
It follows from (23), (24) and (26) that
$$
P^k_{1s}(x_k,\la)=\de_{1s}+O(\rho^{-1}),
\quad \rho\in\Pi_\de,\;|\rho|\to\iy,\; x_k\in(0,T_k].                         \eqno(27)
$$
According to (11) and (25),
$$
P^k_{1s}(x_k,\la)=(-1)^{s-1}\Big(\Big(\vv_k(x_k,\la)\tilde S^{(2-s)}_k(x_k,\la)
-\tilde\vv^{(2-s)}_{k}(x_k,\la)S_k(x_k,\la)\Big)
$$
$$
+(M_k(\la)-\tilde M_k(\la))\vv_k(x_k,\la)\tilde\vv^{(2-s)}_k(x_k,\la)\Big).
$$
It follows from (24) that for each fixed $x_k$, the functions
$P^k_{1s}(x_k,\la)$ are entire in $\la$ of order $1/2.$ Together with
(27) this yields $P^k_{11}(x_k,\la)\equiv 1$, $P^k_{12}(x_k,\la)\equiv 0.$
Substituting these relations into (26) we get $\vv_{k}(x_k,\la)\equiv
\tilde\vv_{k}(x_k,\la)$ for all $x_k$ and $\la,$ and consequently,
$q_k(x_k)=\tilde q_k(x_k)$ a.e. on $[0,T_k],\; h_k=\tilde h_k\,.$
Theorem 2 is proved.

\smallskip
Using the method of spectral mappings [6] for the Sturm-Liouville operator on
the edge $e_k$ one can get a constructive procedure for finding $q_k$ and $h_k$.

\smallskip
Now we study the following auxiliary inverse problem
on the cycle $e_0$, which is called IP(0). Consider the boundary value problem
$B$ of the form (4)-(6), where the parameters of $B_0$ are defined by (3), and
$\al,\be$ are known.

\smallskip
{\bf IP(0).} Given $a(\la), d(\la)$ and $\Omega$, construct
$q(x),\, x\in[0,T], h, \ga_j$ and $\eta_j,\, j=\overline{1,N-1}.$

\smallskip
This inverse problem is a generalization of the classical periodic
inverse problem. Moreover, for the standard matching conditions
($\al_j=\be_j=1, h_j=0$), IP(0) coincides with the classical
periodic inverse problem.

This inverse problem IP(0) was solved in [7], where the following
theorem is established.

\smallskip
{\bf Theorem 3. }{\it The specification $a(\la), d(\la)$ and $\Omega$
uniquely determines $q(x), h, \ga_j$ and $\eta_j,\; j=\overline{1, N-1}.$
The solution of IP(0) can be found by the following algorithm.}

\smallskip
{\bf Algorithm 1. }\\
1) Construct $D(\la)=a(\la)+(1+\al\be).$\\
2) Find zeros $\{z_n\}_{n\ge 1}$ of the entire function $d(\la).$\\
3) Calculate $Q(z_n)$ via
$Q(z_n)=\om_n\sqrt{D^2(z_n)-4\al\be}.$\\
4) Construct $d_1(z_n)$ by
$d_1(z_n)=(D(z_n)+Q(z_n))/(2\alpha).$\\
5) Find $\dot d(z_n).$\\
6) Calculate the Weyl sequence $\{M_n\}_{n\ge 1}$ via
$M_n=-d_1(z_n)/\dot d(z_n).$\\
7) From the given data $\{z_n, M_n\}_{n\ge 1}$ construct $q(x), \ga_j, \eta_j$,
$j=\overline{1,N-1},$ by solving the inverse Dirichlet problem with discontinuities
inside the interval (see [9]).\\
8) Find $S(T,\la),$ $S'(T,\la)$ and $C(T,\la).$\\
9) Calculate $h,$ using (7).

\medskip
{\bf 2.4. } Let us go on to the solution of Inverse problem 1.
Firstly, we give the proof of Theorem 1. Assume that
$\Lambda_k=\tilde\Lambda_k,\; k=\overline{0,r},$
$\Lambda_{{\nu_1,\ldots,\nu_p}}=\tilde\Lambda_{{\nu_1,\ldots,\nu_p}}$,
$p=\overline{2,N},$ $1\le\nu_1<\ldots<\nu_p\le r,$ $\nu_j\in\xi,$
and $\Omega=\tilde\Omega.$ Then one has
$$
\Delta_k(\la)\equiv \tilde\Delta_k(\la),\quad k=\overline{0,r},
$$
$$
\Delta_{\nu_1,\ldots,\nu_p}(\la)\equiv
\tilde\Delta_{\nu_1,\ldots,\nu_p}(\la),\quad
p=\overline{2,N},\;1\le\nu_1<\ldots<\nu_p\le r,\;\nu_j\in\xi.
$$
Moreover, according to (3), $\ga_j=\tilde\ga_j$, $j=\overline{1,N-1},$ and
$\al=\tilde\al,$ $\be=\tilde\be.$
Using Theorem 2, we get $q_k(x_k)=\tilde q_k(x_k)$ a.e. on $[0,T_k]$
and $h_k=\tilde h_k$, $k=\overline{1,r},$ and consequently,
$$
C_k(x_k,\la)\equiv\tilde C_k(x_k,\la),\; S_k(x_k,\la)\equiv\tilde S_k(x_k,\la),
\;\vv_k(x_k,\la)\equiv\tilde\vv_k(x_k,\la),\quad k=\overline{1,r}.             \eqno(28)
$$
By virtue of (15), (17) and (28) one has
$$
\sigma(\la)\equiv\tilde\sigma(\la),\quad
\sigma_{\nu_1,\ldots,\nu_p}(\la)\equiv\tilde\sigma_{\nu_1,\ldots,\nu_p}(\la),
\quad \Omega_j(\la)\equiv\tilde\Omega_j(\la),
\quad \Omega_j^0(\la)\equiv\tilde\Omega_j^0(\la),\quad j=\overline{1,r}.
$$
Using (14) and (16), we obtain, in particular, $a_0(\la)=\tilde a(\la)$,
$a_1(\la)=\tilde a_1(\la).$ In view of (15), this yields
$a(\la)=\tilde a(\la),\; d(\la)=\tilde d(\la).$
It follows from Theorem 3 that $q_k(x_k)=\tilde q_k(x_k)$ a.e. on $[0,T_k],$
$k=\overline{r+1,r+N},$ and $h=\tilde h, \eta_j=\tilde\eta_j,\, j=\overline{1,N-1}.$
Taking (3) into account, we get $H=\tilde H.$ Theorem 1 is proved.

\smallskip
The solution of Inverse problem 1 can be constructed by the
following algorithm.

\smallskip
{\bf Algorithm 2. } Given $\Lambda_k,\, k=\overline{0,r},\,
\Lambda_{{\nu_1,\ldots,\nu_p}},\, p=\overline{2,N},\,
1\le\nu_1<\ldots<\nu_p\le r,\, \nu_j\in\xi$, and $\Omega.$\\
1) Construct $\Delta_k(\la)$ and $\Delta_{\nu_1,\ldots,\nu_p}(\la).$\\
2) Calculate $\ga_j,\; j=\overline{1,N-1},$ $\al$ and $\be,$ using (3).\\
3) For each fixed $k=\overline{1,r},$ solve the inverse problem IP(k)
and find $q_k(x_k),\; x_k\in[0,T_k]$ on the edge $e_k$ and $h_k$.\\
4) For each fixed $k=\overline{1,r},$ construct $C_k(x_k,\la),
S_k(x_k,\la)$ and $\vv_k(x_k,\la),$ $x_k\in[0,T_k].$\\
5) Calculate $a(\la)$ and $d(\la),$ using (14), (15) and (16).\\
6) From the given $a(\la), d(\la)$ and $\Omega,$ construct
$q_k(x_k),\; [0,T_k],$
$k=\overline{r+1,r+N},$ $h$ and $\eta_j$, $j=\overline{1,N-1}.$\\
7) Find $H,$ using (3).

\medskip
{\bf Acknowledgment.} This work was supported by Grant 1.1436.2014K of the Russian 
Ministry of Education and Science and by Grant 13-01-00134 of Russian Foundation for 
Basic Research.

\begin{center}
{\bf REFERENCES}
\end{center}
\begin{enumerate}
\item[{[1]}] Freiling G. and Yurko V.A., Inverse Sturm-Liouville
     Problems and their Applications. NOVA Science Publishers, New York, 2001.
\item[{[2]}] Pokornyi Yu.V. and Borovskikh A.V., Differential equations on networks
     (geometric graphs). J. Math. Sci. (N.Y.) 119, no.6 (2004), 691-718.
\item[{[3]}] Belishev M.I., Boundary spectral inverse problem on a class
     of graphs (trees) by the BC method. Inverse Problems 20 (2004), 647-672.
\item[{[4]}] Yurko V.A., Inverse spectral problems for Sturm-Liouville
     operators on graphs. Inverse Problems 21 (2005), 1075-1086.
\item[{[5]}] Yurko V.A. Uniqueness of recovering differential operators on
     hedgehog-type graphs. Advances in Dynamical Systems and Applications,
     vol.4, no.2 (2009), 231-241.
\item[{[6]}] Yurko V.A., Method of Spectral Mappings in the Inverse
     Problem Theory, Inverse and Ill-posed Problems Series. VSP, Utrecht, 2002.
\item[{[7]}] Yurko V.A. Quasi-periodic boundary value problems with discontinuity
     conditions inside the interval. Schriftenreihe des Fachbereichs Mathematik,
     SM--DU--767, Universit\"{a}t Duisburg-Essen, 2013, 7pp.
\item[{[8]}] M.A. Naimark. Linear differential operators. 2nd ed., Nauka, Moscow, 1969;
     English transl. of 1st ed., Parts I,II, Ungar, New York, 1967, 1968.
\item[{[9]}] Yurko V.A. Integral transforms connected with discontinuous boundary
     value problems, Integral Transforms and Special Functions, vol.10, no.2 (2000), 141-164.
\end{enumerate}

\vspace{0.2cm}

\noindent Yurko, Vjacheslav,\\
Department of Mathematics, Saratov University,\\
Astrakhanskaya 83, Saratov 410026, Russia\\
e-mail: yurkova@info.sgu.ru

\end{document}